\documentstyle[12pt]{article}
\begin{document}
\newcommand{\ced}{\hfill$\Box$}
\def\beq{\begin{eqnarray}}
\def\enq{\end{eqnarray}}
\def\beqq{\begin{eqnarray*}}
\def\enqq{\end{eqnarray*}}
\def\vs{\vspace*{2mm}}
\def\ss{\supseteq}
\def\nss{\not\supseteq}
\def\bs{\subseteq}
\def\nbs{\not\subseteq}
\def\no{\noindent}
\def\nin{\not\in}
\def\ne{\not=}
\def\tra{a contradiction}
\def\la{\lambda_n}
\def\I{\hat{I}}
\def\i{\hat{i}}
\begin{center}
{\bf \large Partition of a Subset  into Two Directed Cycles\vspace*{1mm}

with Partial Degrees

}
\end{center}
\begin{center}
Hong Wang

Department of Mathematics

The University of Idaho

Moscow, Idaho, USA 83844
\end{center}
\parbox[t]{0.5in}{\hfill}
\parbox[t]{5.5in}{\small
Let $D=(V,A)$ be a directed  graph of order $n\geq 6$. Let $W$ be a subset of $V$ with $|W|\geq 6$.  Suppose that every vertex of $W$ has degree at least $(3n-3)/2$ in $D$. Then for any integer partition $|W|=n_1+n_2$ with $n_1\geq 3$ and $n_2\geq 3$, $D$ contains two disjoint directed cycles $C_1$ and $C_2$ such that $|V(C_1)\cap W|=n_1$ and $|V(C_2)\cap W|=n_2$.  We conjecture that for any integer partition $|W|=n_1+n_2+\cdots +n_k$ with $k\geq 3$ and $n_i\geq 3(1\leq i\leq k)$, $D$ contains $k$ disjoint directed cycles $C_1,C_2,\ldots , C_k$ such that $|V(C_i)\cap W|=n_i$ for all $1\leq i\leq k$. The degree condition is sharp in general. }

\section{Introduction}
We discuss only finite simple graphs and strict directed graphs. The terminology and notation concerning graphs is that of \cite{Bondy}, except as indicated. A set of graphs or directed graphs is said to be disjoint if no two of them have any common vertices.  El-Zahar \cite{Zahar} conjectured that if $G$ is a graph of order $n=n_1+\cdots+n_k$ with $n_i\geq 3(1\leq i\leq k)$ and the minimum degree of $G$ is at least $\lceil n_1/2\rceil+\cdots+\lceil n_k/2\rceil$ then $G$ contains $k$ disjoint cycles $C_1,\ldots, C_k$ of orders $n_1,\ldots, n_k$, respectively. This conjecture has been verified for a number of cases (for example, see, \cite{Corradi}, \cite{Zahar} and \cite{Wang2}). Sauer and Spencer in their work \cite{Sauer} conjectured that if the
minimum degree of $G$ is at least $2n/3$ then $G$ contains every
graph of order $n$ with maximum degree of at most 2. This
conjecture was proved by Aigner and Brandt \cite{Aigner}. Shi \cite{Shi}  showed that if $G$ is 2-connected, $W\bs V(G)$ and
$d(x)\geq n/2$ for each $x\in W$ then $G$ contains a cycle passing
through all the vertices of $W$.
This generalizes the classic result of Dirac \cite{Dirac}.\vs

In \cite{Wang4}, we proposed a conjecture  on disjoint cycles of a graph $G$ if we only know a degree condition on vertices in a subset of $V(G)$:\vs
 
\no {\bf Conjecture $A$} \cite{Wang4} {\em Let $G$ be a graph of order $n\geq
3$. Let $W$ be a subset of $V(G)$ with $|W|\geq 3k$ where $k$ is a
positive integer. Suppose that $d(x)\geq 2n/3$ for each $x\in W$.
Then for any integer partition $|W|=n_1+\cdots +n_k$ with $n_i\geq
3(1\leq i\leq k)$, $G$ contains $k$ disjoint cycles
$C_1,\ldots, C_k$ such that $|V(C_i)\cap W|=n_i$ for all $1\leq
i\leq k$.}\vspace*{2mm}

In \cite{Wang4}, we showed that $G$ contains $\lfloor |W|/3\rfloor$ disjoint cycles such that each of them contains at least three vertices of $W$.  This conjecture was verified in \cite{Wang5} for the case $k=2$, i.e., $|W|=n_1+n_2$. This result further supports Conjecture $A$. In \cite{Little} and \cite{Wang3}, we proved the following theorem.\vspace*{2mm}

\no {\bf Theorem $B$ } $($\cite{Little},\cite{Wang3}$)$  {\em Let $D$ be a directed  graph of order $n\geq 4$. Suppose that the  minimum degree of $D$ is at least $(3n-3)/2$. Then for any two integers $s$ and $t$ with $s\geq 2$, $t\geq 2$ and $s+t\leq n$,  $D$ contains two disjoint 
directed cycles of orders $s$ and $t$, respectively.}\vspace*{2mm}

We proposed the following conjecture \cite{Wang1}.\vspace*{2mm}

\no {\bf Conjecture $C$}  \cite{Wang1} {\em Let $D$ be a directed  graph of order $n\geq 4$. Suppose that the  minimum degree of $D$ is at least $(3n-3)/2$. Then for any directed graph $H$ with $\Delta^+(H)\leq 1$ and $\Delta^-(H)\leq 1$, $D$ contains a directed subgraph isomorphic to $H$.}\vspace*{2mm}

To support this conjecture, we proved the following theorem:\vspace*{2mm}

\no {\bf Theorem $D$ } \cite{Wang1} {\em If $D$ is a directed graph of order $n\geq 4$ with minimum degree at least $(3n-3)/2$, then $D$ contains $\lfloor n/3\rfloor$ disjoint directed triangles.}\vspace*{2mm} 

Continuing the work along this line, we propose the following conjecture:\vs

\no {\bf Conjecture $E$} {\em Let $D=(V,A)$ be a directed  graph of order $n\geq 6$. Let $W$ be a subset of $V$.  Suppose that every vertex of $W$ has degree  at least $(3n-3)/2$ in $D$. Then for any integer partition $|W|=n_1+n_2+\cdots +n_k$ with $k\geq 2$ and $n_i\geq 3(1\leq i\leq k)$, $D$ contains $k$ disjoint directed cycles $C_1,C_2,\ldots , C_k$ such that $|V(C_i)\cap W|=n_i$ for all $1\leq i\leq k$. }\vspace*{2mm}

To demonstrate that the degree condition in Conjecture $E$ is sharp in general for each $n\geq 6$, we construct the following directed graph $B_n$ of order $n$. For each positive integer $k$,  we use $K_k^*$
to denote the complete directed graph of order $k$, i.e., $K_k^*$ contains  both $(x,y)$ and $(y,x)$ for any two distinct vertices $x$ and $y$ of $K_k^*$. The directed graph $B_n$ consists of two disjoint complete directed subgraphs $D'$ and $D''$ of order $\lfloor n/2\rfloor$ and $\lceil n/2\rceil$, respectively, and all the arcs $(x,y)$ with $x\in V(D')$ and $y\in V(D'')$. Then the minimum degree of $B_n$ is $\lfloor (3n-4)/2\rfloor$. Clearly, $B_n$ does not contain a directed cycle of order greater than $\lceil n/2\rceil$.\vs

In \cite{Wang6}, we showed that if $|W|\geq 3k$ then $D$ contains $k$ disjoint directed cycles covering $W$ such that each of them contains at least three vertices of $W$. In this paper, we prove the following result:\vspace*{2mm}

\no {\bf Theorem $F$} {\em Let $D=(V,A)$ be a directed  graph of order $n\geq 6$. Let $W$ be a subset of $V$.  Suppose that  every vertex of $W$ has degree  at least $(3n-3)/2$ in $D$.  Then for any integer partition $|W|=n_1+n_2$ with $n_1\geq 3$ and $n_2\geq 3$, $D$ contains two disjoint directed cycles $C_1$ and $C_2$ such that $|V(C_1)\cap W|=n_1$ and $|V(C_2)\cap W|=n_2$. }\vspace*{2mm}

For convenience, we need some special terminology and notation. Let $D$ be a digraph and $G$ a 
graph with the same vertex set, i.e., $V=V(D)=V(G)$.  We use $V(D)$ and $A(D)$ to denote the vertex set and the arc set of
$D$, respectively. Use $E(G)$ to denote the edge set of $G$ and let $e(G)=|E(G)|$. Let  $x\in V$. The out-degree, in-degree and degree of $x$ in $D$ are denoted by $a^+(x)$, $a^-(x)$ and $a(x)$, respectively.  Let each of $X$ and $Y$ be a subset of $V$ or a sequence of some vertices of $V$ or a subgraph of $G$ or a subdigraph of $D$. We use $N_D^+(x, Y)$ denote the set of all the vertices $y$ that are contained in $Y$ with $(x,y)\in A(D)$. We similarly define $N_D^-(x, Y)$ and let
$N_D(x, Y)=N_D^{+}(x, Y)\cup N_D^{-}(x, Y)$. We also define $a^{+}(x, Y)=|N^{+}(x,Y)|$,
$a^{-}(x,Y)=|N^{-}(x,Y)|$ and $a(x,Y)=a^{+}(x,Y)+a^{-}(x,Y)$.  Thus $a^+(x,D)=a^+(x)$, $a^-(x,D)=a^-(x)$ and $a(x,D)=a(x)$.  Define $a(X,Y)=\sum_{u}a(x,Y)$, where $u$ runs over all the vertices of $X$. 

We define $N_G(x,Y)$ to be the set of all the vertices $y$ in $Y$ with $xy\in E$ and let $d(x,Y)=|N_G(x,Y)|$. Hence $d(x,G)=d(x,D)=d(x,V)=d(u)$ which is the degree
of $u$ in $G$. For a subset $U\subseteq V$, $G[U]$ is the subgraph of $G$
induced by $U$ and $D[U]$ is the subdigraph
of $D$ induced by $U$.  Define $d(X,Y)=\sum_{u}d(x,Y)$, where $u$ runs over all the vertices of $X$.  

We define $I(xy,X)=N_G(x,X)\cap N_G(y,X)$ and let $i(xy,X)=|I(xy,X)|$.  We define $I(\vec{xy},X)=N_D^+(x,X)\cap N_D^-(y,X)$ and let $i(\vec{xy},X)=|I(\vec{xy},X)|$.

A graph or a digraph is said to be traceable if it contains a hamiltonian path or  directed hamiltonian path, respectively. 
A graph or  a digraph is called hamiltonian if it
contains a hamiltonian cycle or  directed hamiltonian cycle, respectively. 

If $C$ is a dicycle of $D$ or a cycle of $G$ and $x$ is a vertex of $C$, we use $x^+$ and $x^-$ to denote the successor and the predeccssor of $x$ on $C$, respectively. If $a$ and $b$ are two vertices of $C$, we use $C[a,b]$ to represent the path of $C$ from $a$ to $b$ along the direction of $C$.  We adopt the notation $C(a,b]=C[a,b]-a$, $C[a,b)=C[a,b]-b$ and $C(a,b)=C[a,b]-a-b$.

For two vertices $x$ and $y$ of $D$, we say that $x$ is adjacent to $y$ in $D$ if either $(x,y)$ or $(y,x)$ is an arc of $D$. 
For any two vertices $x$  and $y$ of $G$, we define $\mu (xy)=1$ if $xy$ is an edge of $G$ and $\mu (xy)=0$ otherwise. For any integer $m$, let $\lambda_m=1$ if $m$ is odd and otherwise $\lambda_m=0$.
\section{Lemmas}
\newcommand{\qed}{\hfill\hbox{\rule[-2pt]{4pt}{8pt}}}
Let $D=(V,A)$ be a digraph. Let $G=(V,E)$ be the {\em underlining graph} of $D$, i.e., $V(G)=V$ and $E(G)=\{xy|(x,y)\in A\ {\rm and}\ (y,x)\in A\}$. We will use the following lemmas which are generally well known. 
\newtheorem{lemma}{Lemma}[section]
\newtheorem{coro}[lemma]{Corollary}
\begin{lemma}
Let $P=x_1x_2\ldots x_k$ be a path of $G$ and $\{y, z\}\bs V-V(P)$.  The following two statements hold:

$(a)$ If $d(y,P)>k/2$ then $\{x_i,x_{i+1}\}\subseteq N_G(y,P)$ for some $i\in\{1,2,\ldots ,k-1\}$ unless $k$ is odd, $d(y,P)=(k+1)/2$ and $N_G(y,P)=\{x_1,x_3,x_5,\ldots, x_k\}$;

$(b)$ If $d(y,P)\geq k/2$, then $P+y$ is traceable;

$(c)$ If $d(x_1y,P)\geq k$, then $P+y$ is traceable;

$(d)$ If $d(yz,P)\geq k+2$ then $P+y+z$ is traceable.
\end{lemma}
{\bf Proof.} To show $(a)$, it is easy to see that if $d(y,x_ix_{i+1})\leq 1$ for all $i\in\{1,\ldots, k-1\}$, then $k$ must be odd and $\{x_1,x_k\}\bs N_G(y,P)$ as $d(y,P)>k/2$ and it follows that $N_G(y,P)=\{x_1,x_3,x_5,\ldots, x_k\}$. 

To show $(b)$, it is easy to see that if $\{x_1,x_k\}\cap N_G(y,P)=\emptyset$, then there exists $i\in\{2,\ldots ,k-2\}$ 
such that $\{x_i,x_{i+1}\}\subseteq N_G(y,P)$. Thus the path $P'=x_1\ldots x_iyx_{i+1}
\ldots x_k$ satisfies the requirement. 

To show $(c)$, it is easy to see that if $d(x_k,x_1y)\geq 1$ then $(c)$ holds. So we may assume that $d(x_k,x_1y)=0$. As $d(x_1y,P)\geq k$, it follows that  $\{x_iy,x_{i+1}x_1\}\bs E$ for some $i\in\{1,\ldots, k-2\}$. Therefore $yx_ix_{i-1}\ldots x_1x_{i+1}x_{i+1}\ldots x_k$ is a hamiltonian path of $P+y$. 

To show $(d)$, we see that $\{yx_i,yx_j,zx_{i+1},zx_{j+1}\}\bs E$ for some $\{i,j\}\bs\{1,2,\ldots, k\}$ with  $i<j$, where $x_{k+1}=x_1$. If $j=k$, $zPy$ is a path of $G$ and if $j<k$ then 
\[ x_1\ldots x_iyx_jx_{j-1}\ldots x_{i+1}zx_{j+1}\ldots x_k\] is a hamiltonian path of $P+y+z$.\qed
\begin{lemma}
Let $C$ be a dicycle of $D$ or a cycle of $G$. Let $y\in V-V(C)$. Suppose that $d(y, C)\geq l(C)/2$. Then either $D[V(C)\cup\{y\}]$ is hamiltonian  if $C$ is a cycle of $D$ and $G[V(C)\cup\{y\}]$ is hamiltonian if $C$ is a cycle of $G$, or $l(C)$ is even and there is a list $x_1,\ldots, x_{l(C)}$ of vertices of $C$ along the direction of $C$ such that $N_G(y,C)=\{x_{2i-1}| 1\leq i\leq l(C)/2\}$. 
\end{lemma}
{\bf Proof.}  If the last statement does not hold, then there exists $x\in V(C)$ such that $\{x,x^+\}\subseteq N_G(y)$ and the lemma follows.\qed
\begin{lemma}
Let $P=x_1x_2\ldots x_k$ be a path of $G$ with $k\geq 3$. If either $a^+(x_1,P)+a^+(x_k,P)\geq k$ or  $a^-(x_1,P)+a^-(x_k,P)\geq k$, then $D[V(P)]$ is hamiltonian. Therefore if $a(x_1,P)+a(x_k,P)\geq 2k-1$ then $D[V(P)]$ is hamiltonian.
\end{lemma}
{\bf Proof.} If $a^+(x_1,P)+a^+(x_k,P)\geq k$, then $\{(x_1,x_{i+1}),(x_k,x_i)\}\bs A$ for some $i\in\{1,\ldots, k-1\}$ and so $(x_1,x_{i+1}, x_{i+2},\ldots, x_k,x_i,x_{i-1},\ldots, x_1)$ is a dicycle of $D[V(P)]$. Similarly, we can easily see that  the lemma holds if $a^-(x_1,P)+a^-(x_k,P)\geq k$. Since $a(x_1,P)+a(x_k,P)\geq 2k-1$ implies either $a^+(x_1,P)+a^+(x_k,P)\geq k$ or  $a^-(x_1,P)+a^-(x_k,P)\geq k$, the lemma follows.\qed
\begin{lemma}
{\rm \cite{Erdos}} Let $P=x_tx_{t-1}\ldots x_1$ be a longest path starting at $x_t$ in $G$. Let $c>t/2$.
Suppose that for each $v\in V(G)$, if there
exists a longest path starting at $x_t$ in $G$ such that the path
ends at $v$ then $d_G(v)\geq c$. Then there exists $r\geq c+1$ such that $N(x_i)\bs
\{x_1,x_2,\ldots, x_r\}$, $G[P]$ has an $x_t$-$x_i$
hamiltonian path containing $x_tx_{t-1}\ldots x_r$ and $d_G(x_i)\geq c$ for all $i\in\{1,2,\ldots,
r-1\}$. Moreover, if $t>r$ then $x_r$ is a cut-vertex of $G$.
\end{lemma}
\section{Proof of Theorem $F$}
Let $D=(V,A)$ be a directed graph of order $n$. Let $W\bs V$ with $|W|\geq 6$. Suppose that $a(x)\geq (3n-2-\lambda_n)/2$ for all $x\in W$. Suppose, for a contradiction, that for some integer partition $|W|=n_1+n_2$ with $n_1\geq 3$ and $n_2\geq 3$, $D$ does not contain two disjoint dicycles $C_1$ and $C_2$ such that $|V(C_1)\cap W|=n_1$ and $|V(C_2)\cap W|=n_2$.  For the proof, we may assume that $D-W$ has no arcs. Let $G=(V, E)$ be the underlining graph of $D$.  

For any subgraph $G'$ of $G$,  we use $W(G')$ to denote $W\cap V(G')$.  For convenince, let $A'=\{xy|\ \{x,y\}\bs  V\ {\rm and\ either}\ (x,y)\in A\ {\rm or}\ (y,x)\in A\}$, i.e.,  $x$ and $y$ are adjacent in $D$ if and only if $xy\in A'$.  For a dicycle or a cycle $C$, we use $l(C)$ to denote the length of $C$ and use $l_w(C)$ to denote $|V(C)\cap W|$ which is called the $W$-{\em length} of $C$. Similarly, we define $l(P)$ and $l_w(P)$ for a dipath or a path $P$. 

For each $x\in W$, we have
\beq
d(x,G)&\geq& (3n-2-\la)/2-|N_D(x)|\\
&\geq&(3n-2-\la)/2-(n-1)=(n-\la)/2.
\enq

For each partition $(X,Y)$ of $V$ and $u\in X$, if $a(u,X)\leq |X|-1$, then we have
\beq
a(u,Y)&\geq& (3n-2-\la)/2-a(u,X)\geq  (3|Y|+|X|-\la)/2;\\
|Y|&\geq& d(u, Y)=a(u,Y)-|N_D(u,Y)|\\
& \geq&(3|Y|+|X|-\lambda_n)/2-|N_D(u, Y)|\\
&\geq&(3|Y|+|X|-\lambda_n)/2-|Y|=(n-\lambda_n)/2.
\enq
 
\begin{lemma}
There is a path $P$ in $G$ such that $V(P)\ss W$.
\end{lemma}
{\bf Proof.} We choose a path $P$ in $G$ with $l_w(P)$ maximal. Subject to this, we 
choose $P$ with $l(P)$ minimal. Assume there exists $x_0\in W-V(P)$.  Say  $P=x_1x_2\ldots x_t$. Then $x_0x_1\nin E$, $I(x_1x_0,V-V(P))=\emptyset$ and so $d(x_1x_0,V-V(P))\leq n-t-1$. By (2), $d(x_1x_0,P)\geq (n-\lambda_n)-(n-t-1)\geq t$. By Lemma 2.1$(c)$, $P+x_0$ is traceable, a contradiction. \qed\vs

By Lemma 3.1, $G$ has two disjoint paths $P_1$ and $P_2$ with $l_w(P_1)=n_1$ and $l_w(P_2)=n_2$. Let $G_1=G[V(P_1)]$ and $G_2=G[V(P_2)]$. We choose $P_1$ and $P_2$ such that
\beq
l(P_1)+l(P_2)\ {\rm is\  minimal.}
\enq
Subject (7), we choose $P_1$ and $P_2$ such that 
\beq
e(G_1)+e(G_2)\ {\rm is\  maximal.}
\enq
Say $P_1=x_1x_2\ldots x_s$, $P_2=y_1y_2\ldots y_t$, $R=V-V(G_1\cup G_2)$ and $r=|R|=n-s-t$. Let $H=G[V(G_1\cup G_2)]$, $D_1=D[V(G_1)]$ and $D_2=D[V(G_2)]$. By (8), we obtain the following lemma.
\begin{lemma}
Let $x$ be an endvertex of a hamiltonian
path of $G_1$ and $y$ an endvertex of a hamiltonian path of $G_2$. If
both $G_1-x+y$ and $G_2-y+x$ are traceable, then
\begin{eqnarray}
d(x, G_1)+d(y,G_2)\geq d(x,G_2)+d(y,G_1)-2\mu (xy).
\end{eqnarray}
In particular, $(9)$ holds if $d(x,G_2)\geq (t+1)/2$ and 
$d(y,G_1)\geq (s+1)/2$ hold.
\end{lemma}
{\bf Proof.} We have 
\begin{eqnarray*}
& &e(G_1-x+y)+e(G_2-y+x)\\
&=& e(G_1)+e(G_2)-d(x,G_1)-d(y,G_2)\\
&+&
d(x,G_2)+d(y,G_1)-2\mu (xy).
\end{eqnarray*}
If both $G_1-x+y$ and $G_2-y+x$ are traceable, then $e(G_1)+e(G_2)\geq e(G_1-x+y)+e(G_2-y+x)$ by (8) and it follows that (9) holds. Note that by Lemma 2.1$(b)$, if $d(x,G_2)\geq (t+1)/2$ and $d(y,G_1)\geq (s+1)/2$, both $G_1-x+y$ and $G_2-y+x$ are traceable. \qed\vs

Similarly, by (8), one can easily see that  the following Lemma 3.3 holds.
\begin{lemma}
Let $x$ and $y$ be the two endvertices of a hamiltonian path of $G_1$ and $u$ and $v$ be the two endvertices of a hamiltonian path of $G_2$.  Suppose that
both $G_1-x-y+u+v$ and $G_2-u-v+x+y$ are taceable. Then 
\beq
d(xy,G_1)+d(uv,G_2)\geq d(xy,G_2)+d(uv,G_1)-2d(xy,uv).
\enq
\end{lemma}

\begin{lemma}
For some $i\in\{1,2\}$, $D_i$ does not have a dicycle of  $W$-length $n_i$ and if $x$ and $y$ are the two endvertices of a hamiltonian path of $G_i$, then $I(\vec{xy},R)=\emptyset$,  $I(\vec{yx},R)=\emptyset$ and so $a(xy,R)\leq 2r$. 
\end{lemma}
{\bf Proof.} If $D_1$ contains a dicycle of $W$-length $n_1$, then $D[V(D_2)\cup R]$ does not contain a dicycle of $W$-length $n_2$, and this implies that the lemma holds with $i=2$. Therefore we may assume that $D_1$ does not contain a dicycle of $W$-length $n_1$. Similarly, we may assume that $D_2$ does not contain a dicycle of $W$-length $n_2$. In particular, $x_1x_s\nin A'$ and $y_1y_t\nin A'$. On the contrary, say the lemma fails. Say w.l.o.g. that $I(\vec{x_1x_s},R)\not=\emptyset$ and $I(\vec{y_1y_t},R)\not=\emptyset$. Since $D$ does not contain the two required dicycles, we see that for some $u\in R$, $I(\vec{x_1x_s},R)=I(\vec{y_1y_t},R)=\{u\}$.  Moreover, $I(\vec{x_sx_1},R)\bs\{u\}$ and $I(\vec{y_ty_1},R)\bs \{u\}$.  This implies that $a(x_1x_s,R)\leq 2r+2$ and $a(y_1y_t, R)\leq 2r+2$. Then for each $vw\in\{x_1x_s,y_1y_t\}$, $a(vw,H)\geq 3n-2-\lambda_n-2r-2=3s+3t+r-4-\lambda_n$. Therefore for each $vw\in\{x_1x_s,y_1y_t\}$, we have 
\beq
d(vw ,H)&\geq& 3s+3t+r-4-\lambda_n-|N_D(v,H)|-|N_D(w,H)|\\
&\geq &3s+3t+r-4-\lambda_n-2(s+t-2)=s+t+r-\lambda_n.
\enq
Since both $G_1$ and $G_2$ are not hamiltonian, $d(x_1x_s,G_1)\leq s-1$ and $d(y_1y_t,G_2)\leq t-1$. Therefore $d(x_1x_s, G_2)\geq s+t+r-\la-(s-1)\geq t+1+r-\la$ and $d(y_1y_t,G_1)\geq  s+t+r-\lambda_n-(t-1)\geq s+1+r-\la$. Let $\sigma=d(x_1x_s,G_2)+d(y_1y_t,G_1)-d(x_1x_s,G_1)-d(y_1y_t,G_2)$. Then $\sigma\geq 4+2(r-\lambda_n)$.  

First, assume that there exist two independent edges between $\{x_1,x_s\}$ and $\{y_1,y_t\}$ in $G$.  Say w.l.o.g. $\{x_1y_t,x_sy_1\}\bs E$. By Lemma 3.2,  $d(x_1,G_1)+d(y_1,G_2)\geq d(x_1,G_2)+d(y_1,G_1)-2\mu(x_1y_1)$ and  $d(x_s,G_1)+d(y_t,G_2)\geq d(x_s,G_2)+d(y_t,G_1)-2\mu(x_sy_t)$.  Thus $2\mu(x_1y_1)+2\mu(x_sy_t)\geq \sigma$. Since $4\geq 2\mu(x_1y_1)+2\mu(x_sy_t)$ and $\sigma\geq 4+2(r-\lambda_n)$, it follows that equality hold in (11) and (12) with $r=\lambda_n=1$.  In particular, $|N_D(z,H)|=s+t-2$ for all $z\in\{x_1,x_s,y_1,y_t\}$.  Hence  $x_1y_2\in A'$ and $y_1x_2\in A'$.  Consequently, $D_1-x_1+y_1$ and $D_2-y_1+x_1$ are hamiltonian, \tra. 

Therefore $G$ does not contain two independent edges between $\{x_1,x_s\}$ and $\{y_1,y_t\}$ and so $d(x_1x_s,y_1y_t)\leq 2$.  Assume for the moment that $d(x_1x_s,G_2)\geq (t+2)$ and $d(y_1y_t,G_1)\geq s+2$. Then $d(x_1x_s,G_2-y_1-y_t)\geq (t-2)+2$ and $d(y_1y_t,G_1-x_1-x_s)\geq (s-2)+2$. By Lemma 2.1$(d)$ and Lemma 3.3, we shall have that $d(x_1x_s,G_1)+d(y_1y_t,G_2)\geq d(x_1x_s,G_2)+d(y_1y_t,G_1)-2d(x_1x_s,y_1y_t)$. Thus $s-1+t-1\geq s+2+t+2-4$, \tra. Therefore we may assume w.l.o.g.that $d(x_1x_s,G_2)\leq t+1$. By (12), $d(x_1x_s,G_2)\geq s+t+r-\lambda_n-(s-1)= t+1+r-\lambda_n$. It follows that $d(x_1x_s,G_1)=s-1$, $d(x_1x_s,G_2)=t+1$, $r=1$, $\lambda_n=1$ (i.e., $n$ is odd), and equality hold in (11) and (12) with $\{v,w\}=\{x_1,x_s\}$ and $N_D(x_1)=N_D(x_s)=V-\{x_1,x_s\}$.  Moreover, $a(x_1x_s,R)=a(x_1x_s,u)=4$, i.e.,  $d(x_1x_s,R)=d(x_1x_s,u)=2$.  As $x_1x_s\nin A'$ and by (1), we see that $d(x_1,G)\geq (n+1)/2$ and $d(x_s,G)\geq (n+1)/2$. As $x_1x_{s-1}\in A'$, $D_1-x_s$ is hamiltonian and so $d(x_s,G_1)\leq (s-1)/2$ by Lemma 2.2 because $D_1$ is not hamiltonian. Similarly, $d(x_1,G_1)\leq (s-1)/2$. It follows that $s$ and $t$ are odd and $d(x_i,G_1)=(s-1)/2$ and $d(x_i,G_2)=(t+1)/2$ for $i\in\{1,s\}$. Suppose that $d(y_1y_t,G_1)\leq s+1$. Then this argument allows us to see that $d(y_j,G_2)=(t-1)/2$ and $d(y_j,G_1)=(s+1)/2$ for $j\in\{1,t\}$. For each $i\in\{1,s\}$ and $j\in\{1,t\}$, applying Lemma 3.2 with $x_i$ and $y_j$ in place of $x$ and $y$, we see that $\mu (x_iy_j)=1$. Thus $d(x_1x_s,y_1y_t)=4$, \tra. 

Therefore $d(y_1y_t,G_1)\geq s+2$. Say w.l.o.g. $d(y_1,G_1)\geq \lceil (s+2)/2\rceil=(s+3)/2$. Applying Lemma 3.2 with $x_1$ and $y_1$ in place of $x$ and $y$, we see that $d(y_1,G_2)\geq (t+1)/2$. Thus $d(y_t, G_2)\leq t-1-(t+1)/2=(t-3)/2$. By (1), $d(y_t,G)\geq (n+1)/2$.  Consequently, $d(y_t,G_1)\geq (n+1)/2-(t-3)/2-\mu (y_tu)\geq (s+3)/2$. Applying Lemma 3.2 with $x_1$ and $y_t$ in place of $x$ and $y$, we obtain a contradiction to (9).\qed\vspace*{2mm}

We now may assume that Lemma 3.4 holds for $i=1$. As $D_1$ is not hamiltonian, $a(x_1x_s,P_1)\leq 2(s-1)$ by Lemma 2.3 and so $a(x_1x_s,G_1+R)\leq 2(s-1)+2r$. Thus 
\beq
& &a(x_1x_s,G_2)\geq 3n-2-\la-a(x_1x_s, G_1+R)\geq 3t+s+r-\la;\nonumber\\
& &d(x_1x_s, G_2)\geq  3t+s+r-\la-|N_D(x_1,G_2)|-|N_D(x_s,G_2)|\geq t+s+r-\la.
\enq
 We assume w.l.o.g. that $a(x_1,G_1+R)\leq a(x_s,G_1+R)$. Then $a(x_1,G_1+R)\leq s+r-1$ and $a(x_1, G_2)\geq (3n-2-\la)/2-(s+r-1)=(3t+s+r-\la)/2$. Thus 
\beq
d(x_1,G_2)&\geq& (3t+s+r-\la)/2-|N_D(x_1, G_2)|\\
&\geq& (3t+s+r-\la)/2-t=(t+s+r-\la)/2\\
&=&(n-\la)/2.
\enq
It follows that $t\geq (n-\la)/2$. Clearly, $(n-\la)/2-1>(t-1)/2$. By Lemma 2.1$(b)$ and Lemma 2.2, we obtain the following corollary:

\begin{coro} For each $v\in V(G_2)$, if $G_2-v$ is traceable then $G_2-v+x_1$ is traceable and if $D[V(G_2-v)]$ is hamiltonian then $D[V(G_2-v+x_1)]$ is hamiltonian.\qed
\end{coro}

\begin{lemma}
$G_2$ is not hamiltonian and $G_2+z$ is not hamiltonian for each $z\in R$. 
\end{lemma}
{\bf Proof.} On the contrary, say the lemma fails. If $G_2$ is a hamiltonian, let $C$ be a hamiltonian cycle and otherwise let $C$ be a hamiltonian cycle of $G_2+z_0$ for some $z_0\in R$. Let $b_1,b_2,\ldots, b_{n_2}$ be a list of vertices in $W(C)$ along the direction of $C$.  Let $h_i=C[b_i,b_{i+1})$ for each $i\in\{1,\ldots, n_2\}$, where the operation in the subscripts is taken in $\{1,2,\ldots, n_2\}$ by modulo $n_2$.  Note that as $D-W$ has no arcs, each $h_i$ has at most two vertices. Let $R'=R-\{z_0\}$ if $l(C)=t+1$ and otherwise let $R'=R$. Say $r'=|R'|$. 

Let $x^*=x_2$ if $x_2\in W$ and otherwise let $x^*=x_3$ with $x_3\in W$. Suppose that $a(x^*x_s,h_i)\leq 2|V(h_i)|$ for all $i\in\{1,2, \ldots, n_2\}$. Then $a(x^*x_s,C)\leq 2l(C)$. Thus $a(x^*x_s,D-V(C))\geq 3n-2-\la-2l(C)= 3s+3r'+l(C)-2-\la$. As $x_1x_s\nin A'$, we obtain that $2(s+r'-1)+2(s+r'-2)\geq a(x^*x_s,D-V(C))\geq 3s+3r'+l(C)-2-\la$ and so $s+r'\geq l(C)+4-\la$. This implies that $n\geq 2l(C)+4-\la$, i.e., $l(C)\leq (n-4+\la)/2$. But $l(C)\geq t\geq (n-\la)/2$, \tra. 

Therefore $a(x^*x_s,h_i)\geq 2|V(h_i)|+1$ for some $i\in\{1,\ldots, n_2\}$. Then either $a^+(x^*,h_i)+a^-(x_s,h_i)\geq |V(h_i)|+1$ or $a^-(x^*,h_i)+a^+(x_s,h_i)\geq |V(h_i)|+1$.  As $|V(h_i)|\leq 2$, it follows that $D[V(P_1[x^*,x_s])\cup V(h_i)]$ has a hamiltonian dicycle $C_1$ of $W$-length $n_1$. Therefore $D-V(C_1)$ does not contain a dicycle  of $W$-length $n_2$.  In particular, $d(x_1,b_{i-1}b_{i+1})\leq 1$ and $a(x_1,b_{i-1}b_{i+1})\leq 2$. Recall that $a(x_1,G_1+R)\leq s+r-1$. By (3) to (6) with $Y=V(G_2)-\{b_{i-1},b_{i+1}\}$ and $x_1$ in place of $u$, we obtain 
\beq
t-2\geq d(x_1,G_2-\{b_{i-1},b_{i+1}\})\geq (n-\la)/2\  {\rm and\ so}\ t\geq (n+4-\la)/2.
\enq
 By (16), we now have
\beq
d(x_1,G_2-V(h_i))\geq (t+s+r-\la)/2-|V(h_i)|.
\enq  
Note that if $C$ is a hamiltonian cycle of $G_2$, then by (7), $V(G_2)\bs W$, i.e., $n_2=t$.  We claim
\beq
{\rm Either}\  a(b_{i-1}b_{i+1}, C-b_i)\geq 2(l(C)-1)\ {\rm or}\ a(b_{i-1}b_{i+1}, R')\geq 2r'+1.
\enq

To see (19), say, for a contradiction, that  $a(b_{i-1}b_{i+1}, R'\cup V(C)-\{b_i\})\leq 2r'+2(l(C)-1)-1$. As $a(x_1,b_{i-1}b_{i+1})\leq 2$, $a(b_{i-1}b_{i+1}, R'\cup V(C)\cup\{x_1\}-\{b_i\})\leq 2r'+2l(C)-1=2(t+r-1)+1$. This implies that $a(u,  R'\cup V(C)\cup\{x_1\}-\{b_i\})\leq t+r-1$ for some $u\in\{b_{i-1},b_{i+1}\}$. By (3) to (6) with $Y=V(G_1)\cup\{b_i\}-\{x_1\}$, we obtain $s\geq (n-\la)/2$. By (17), $s+t\geq  (n-\la)/2+ (n+4-\la)/2\geq n+1$, \tra.  Hence (19) holds. 

We break into the following four cases. Cases 2 to 4 are similar to Case 1 with more subtle details. \vs

 \noindent Case 1. $b_i^+=b_{i+1}$ and $b_i^-=b_{i-1}$.

 Then $h_i=b_i$. Let $P=C-b_i$.  If the latter of (19) holds, $a(z', b_{i-1}b_{i+1})\geq 3$ for some $z'\in R'$ and we let $C'$ be a hamiltonian dicycle of $D[V(P)\cup\{z'\}]$. Otherwise the former holds and we let $C'$ be a hamiltonian dicycle of $D[V(P)]$ by Lemma 2.3. By (18), $d(x_1,C')\geq (t+s+r-\la)/2-1=(t+s-2+r-\la)/2>l(C')/2$. Thus $D[V(C')\cup\{x_1\}]$ is hamiltonian by Lemma 2.2, \tra. \vs

\noindent Case 2. $b_i^+=b_{i+1}$ and $b_i^{--}=b_{i-1}$.

Then $h_i=b_i$. In this case, as $b_i^-\nin W$, $C$ is not a hamiltonian cycle of $G_2$ by (7). Thus $l(C)=t+1$ with $V(C)\cap R=\{z_0\}$ and so $r\geq 1$ . As in Case 1, we see that we may assume that $D[V(C)-\{b_i\}]$ is not hamiltonian. Thus $b_i^-b_{i+1}\nin A'$.  If the former of (19) holds, then $a(b_{i-1}b_{i+1}, C-b_i-b_i^-)\geq 2(l(C)-2)$ and we let $C''$ be a hamiltonian dicycle of $D[V(C-b_i-b_i^-)]$ by Lemma 2.3. Otherwise $D[V(C-b_i-b_i^-)\cup\{z''\}]$ is hamiltonian for some $z''\in R'$ with $a(b_{i-1}b_{i+1}, z'')\geq 3$ and we let $C''$ be one of its hamiltonian dicycles.  If $x_1b_i^-\nin E$, then $d(x_1,C'')\geq (t+s+r-\la)/2-1\geq (l(C'')+1)/2$ and so $D[V(C'')\cup\{x_1\}]$ is hamiltonian, a contradiction. Therefore $x_1b_i^-\in E$. Then $x_1b_{i+1}\nin A'$ for otherwise $D[V(C)\cup\{x_1\}-\{b_i\}]$ is hamiltonian. By (14), we see that $d(x_1,C'')\geq d(x_1,G_2)-2\geq (3t+s+r-\la)/2-(t-1)-2= (t+s-2+r-\la)/2\geq (l(C'')+1)/2$.   Again, $D[V(C'')\cup\{x_1\}]$ is hamiltonian,  \tra.\vs

\no Case 3. $b_i^{++}=b_{i+1}$ and $b_i^-=b_{i-1}$.

Then $h_i=b_ib_i^+$. In this case, $l(C)=t+1$ with $V(C)\cap R=\{z_0\}$.  Let $P'=C-b_i-b_i^+$. First, assume that the latter of (19) holds. Then $D[V(P')\cup\{u\}]$ has a hamiltonian dicycle $L$ for some $u\in R'$ with $a(b_{i-1}b_{i+1},u)\geq 3$. Thus $r'\geq 1$ and so $r\geq 2$. Clearly, $l(L)\leq t$ and $d(x_1,L-u-z_0)=d(x_1,G_2-b_i-b_i^+)\geq (t+s+r-\la)/2-2=(t+s+r-4-\la)/2\geq l(L)/2$.  As $D[V(L)\cup\{x_1\}]$ is not hamiltonian and by Lemma 2.2, it follows  that $l(L)=t$,  $d(x_1,G_2-b_i-b_i^+)=d(x_1,L)=t/2$, $d(x_1, N_L(z_0))=2$ and so $D[V(L)\cup\{x_1\}-\{z_0\}]$ is hamiltonian, \tra. Hence the latter of (19) does not hold. Note that this argument allows us to see that $D[V(P')]$ does not have a hamiltonian dicycle in both Case 3 and Case 4. 

Therefore  $a(b_{i-1}b_{i+1}, R')\leq 2r'$, and by Lemma 2.3, $a(b_{i-1}b_{i+1},P')\leq 2(|V(P')|-1)$. As $a(x_1, b_{i-1}b_{i+1})\leq 2$,  it follows that $a(b_{i-1}b_{i+1}, R'\cup V(C)\cup\{x_1\}-\{b_i,b_i^+\})\leq 2(t+r-2)$. Thus $a(v, R'\cup V(C)\cup\{x_1\}-\{b_i,b_i^+\})\leq t+r-2$ for some $v\in\{b_{i-1},b_{i+1}\}$. By (3) to (6) with $Y=V(G_1)\cup\{b_i,b_i^+\}-\{x_1\}$, we obtain $s+1\geq (n-\la)/2$. By (17), $n\geq s+t+1\geq (n-\la)/2+(n+4-\la)/2\geq n+1$, \tra.  \vs

\no Case 4. $b_i^{++}=b_{i+1}$ and $b_i^{--}=b_{i-1}$. 

In this case, $h_i=b_ib_i^+$, $l(C)=t+1$ and $r\geq 1$.  As noted in Case 3,  we may assume that $D[V(P')]$ is not hamiltonian for otherwise $D[V(P')\cup\{x_1\}]$ is hamiltonian, where $P'= C-b_i-b_i^+$.   Similarly, $D[V(C)-\{b_i,b_i^-\}]$ is not hamitonian. Thus $b_{i+1}b_i^-\nin A$ and $b_{i-1}b_i^+\nin A'$. Set $P''=C-b_i^- -b_i-b_i^+$. 

If the former of (19) holds, then $a(b_{i-1}b_{i+1}, P'')\geq 2(l(C)-1)-4=2|V(P'')|$ and by Lemma 2.3, we let $Q$ be a hamiltonian dicycle of $D[V(P'')]$. Otherwise $D[V(P'')\cup\{c\}]$ is hamiltonian for some $c\in R'$ with $a(c, b_{i-1}b_{i+1})\geq 3$, and we let $Q$ be a hamiltonian dicycle of $D[V(P'')\cup\{c\}]$.  If $x_1b_i^-\nin E$,  $d(x_1,Q)\geq (t+s+r-\la)/2-2=(t+s-4+r-\la)/2>l(Q)/2$ and so $D[V(Q+x_1)]$ is hamiltonian by Lemma 2.2, \tra. Therefore $x_1b_i^-\in E$. In this situation, $x_1b_{i+1}\nin A'$. By (14), $d(x_1,G_2)\geq (t+s+r+2-\la)/2$.  Thus $d(x_1,Q)\geq  (t+s+r+2-\la)/2-3=(t+s-4+r-\la)/2>l(Q)/2$ and so $D[V(Q+x_1)]$ is hamiltonian, \tra.\qed\vs

By Lemma 3.6, we obtain
\beq
& &d(y_1y_t, G_2)\leq t-1, I(y_1y_t,R)=\emptyset\ {\rm and}\ d(y_1y_t, G_2+R)\leq r+t-1.
\enq
Since we have $d(y_1y_t,G)\geq 3n-2-\la-|N_D(y_1)|-|N_D(y_t)|$, it follows that
\beq
d(y_1y_t,G_1)&\geq& 3n-2-\la-|N_D(y_1)|-|N_D(y_t)|-d(y_1y_t,G_2+R)\\
&\geq & 3s+2t+2r-1-\la-|N_D(y_1)|-|N_D(y_t)|\\
&\geq &3s+2t+2r-1-\la-2(s+t+r-1)\\
&=&s+1-\la\geq s.
\enq
From (21) and (22), it follows that
\beq
& &{\rm If}\ d(y_1y_t,G_1)=s,\  {\rm then}\ \la=1\ {\rm and}\   |N_D(y_1)|=|N_D(y_t)|=n-1;\\
& &{\rm If}\ d(y_1y_t,G_1)=s+1,\  {\rm then}\ y_1y_t\in A'\ {\rm and}\   |N_D(y_1)|+|N_D(y_t)|\geq 2n-3.
\enq
 We claim
\beq
d(x_1x_s,y_1y_t)\leq 2, d(y_1y_t,G_1)\leq s+1\ {\rm and}\ y_1y_t\in A'.
\enq

{\em Proof of $(27)$.} On the contray, say (27) fails. Since $G_1$ is not  hamiltonian, $d(x_1x_s,G_1)\leq s-1$. First, assume that $d(x_1x_s,y_1y_t)\geq 3$. Then there exist two independent edges of $G$ between $\{x_1,x_s\}$ and $\{y_1,y_t\}$. Say w.l.o.g. $\{x_1y_t,x_sy_1\}\bs E$.  If $y_ty_2\in A'$ then $D[V(G_2)\cup\{x_1\}-\{y_1\}]$ is hamiltonian by Corollary 3.5 and so $y_1x_2\nin A'$. That is, either $y_ty_2\nin A'$ or $y_1x_2\nin A'$. By (22), we obtain  $d(y_1y_t,G_1)\geq s+2-\la$. With (13) and (20), we obtain
\beqq
& &d(x_1x_s,G_2)+d(y_1y_t,G_1)-d(x_1x_s,G_1)-d(y_1y_t,G_2)\\
& &\geq (t+s+r-\la)+(s+2-\la)-(s-1)-(t-1)\\
& &=s+4+r-2\la\geq 5.
\enqq
Thus $d(x,G_2)+d(y,G_1)-d(x,G_1)-d(y,G_2)\geq 3$ for some $(x,y)\in\{(x_1,y_1),(x_s,y_t)\}$. This contradicts Lemma 3.2. Hence $d(x_1x_s,y_1y_t)\leq 2$. Thus by (13), $d(x_1x_s,G_2-y_1-y_t)\geq t+s+r-\la-2\geq (t-2)+2$. Next, assume $d(y_1y_t,G_1)\geq s+2$.  Then $d(y_1y_t,G_1-x_1-x_s)\geq s+2-2=(s-2)+2$. By Lemma 2.1$(d)$, both $G_1-x_1-x_s+y_1+y_t$ and $G_2-y_1-y_t+x_1+x_s$ are traceable. By Lemma 3.3, we would have that $(s-1)+(t-1)\geq (t+s+r-\la)+(s+2)-4$, which is impossible. Hence $d(y_1y_t,G_1)\leq s+1$. Together with (24), (25) and (26), we obtain $y_1y_t\in A'$. Hence (27) holds. \ced\vs

We now choose  $P_2=y_1y_2\ldots y_t$ in $G_2$ such that
\beq
 d(y_1y_t,G_2)\ {\rm is\  minimal}.
\enq
 
\begin{lemma}
For each $v\in\{y_1,y_t\}$, if $d(v,G_2)\leq (t-1)/2$ then $d(v,G_1)\leq (s-1)/2$. 
\end{lemma}
{\bf Proof.} On the contrary, say $d(y_1,G_2)\leq (t-1)/2$ and $d(y_1, G_1)\geq s/2$. By Corollary 3.5, $G_2-y_1+x_1$ is traceable. First, assume that $G_1-x_1+y_1$ is tracable. Then by Lemma 3.2, $(s-2)+(t-1)/2\geq (t+s+r-\la)/2+\lceil s/2\rceil-2\mu (x_1y_1)$.  This implies that $r=0$, $\la=1$, $d(x_1,G_1)=s-2$, $d(y_1,G_2)=(t-1)/2$,  $d(y_1,G_1)=s/2$ and $x_1y_1\in E$. Thus $s$ is even, $s\geq 4$ and $N_G(x_1,G_1)=\{x_2,\ldots, x_{s-1}\}$. It follows that $d(y_t,G_2)\leq (t-1)/2$ by (20),  and by (24), $d(y_t,G_1)\geq s/2$. Clearly,  $d(x_s,G_1)=1$ as $d(x_1x_s,G_1)\leq s-1$, and $G_1-x_s+y_t$ is traceable.  As $x_1x_s\nin A'$ and by (1), $d(x_s,G)\geq (t+s+1)/2$ and so $d(x_s,G_2)\geq (t+s-1)/2$. Thus $G_2-y_1+x_s$ is traceable.  Applying Lemma 3.2 with $x_s$ and $y_1$ in place of $x$ and $y$, we obtain $1+(t-1)/2\geq (t+s-1)/2+s/2-2\mu(x_sy_1)$, which is impossible.  Hence $G_1-x_1+y_1$ is not traceable.  Therefore $d(y_1,x_2x_s)=0$ and $d(y_1,P_1-\{x_1,x_2,x_s\})\leq \lfloor (s-2)/2\rfloor$ by Lemma 2.1$(a)$.   As $d(y_1,G_1)\geq s/2$, it follows that  $d(y_1,G_1)=s/2$ with $N_G(y_1,G_1)=\{x_1,x_3,x_5,\ldots, x_{s-1}\}$.  Then $N_G(x_2,G_1-x_1)\bs\{x_3,x_5,\ldots, x_{s-1}\}$ and $N_G(x_s,G_1)\bs\{x_3,x_5,\ldots, x_{s-1}\}$ for otherwise $G_1-x_1+y_1$ is traceable. Thus $d(x_2,G_1)\leq s/2$ and $d(x_s,G_2)\leq (s-2)/2$. We claim
\begin{eqnarray}
d(x_s,G_2)\leq (t-2)/2\ {\rm and\ if}\ x_2\in W\ {\rm then}\ d(x_2,G_2)\leq (t-1)/2.
\end{eqnarray}

To see (29), we have that $d(x_2,G_1)\leq s/2$, $d(y_1,x_2x_s)=0$ and $G_1-x_i+y_1$ is traceable for each $i\in\{2,s\}$. If $d(x_s,G_2)\geq (t-1)/2$, then $G_2-y_1+x_s$ is traceable by Lemma 2.1$(b)$ and so $(s-2)/2+(t-1)/2\geq (t-1)/2+s/2$ by Lemma 3.2, which is impossible. Hence $d(x_s,G_2)\leq (t-2)/2$.  If $x_2\in W$ and  $d(x_2,G_2)\geq t/2$, then $G_2-y_1+x_2$ is traceable by Lemma 2.1$(b)$ and so $s/2+(t-1)/2\geq t/2+s/2$ by Lemma 3.2, which is impossible. Therefore (29) holds.

Then $d(x_s,R)\geq (n-\la)/2-(s-2)/2-(t-2)/2=(r+4-\la)/2$. We also have that $d(y_1,R)\geq (n-\la)/2-s/2-(t-1)/2=(r+1-\la)/2$. We obtain that $d(x_sy_1,R)\geq r+5/2-\la$ and so $i(x_sy_1,R)\geq 2$. Note that $r\geq 2$.  If $x_2\in W$ and $a(x_2,N_G(y_1,R))=0$ then $d(x_2,G)\geq (n+4-\la)/2$ by (1) and so $i(x_2y_1,R)\geq 2$, \tra.  Therefore if $x_2\in W$ then $a(x_2, N_G(y_1,R))>0$.  It follows that $D_1-x_1+y_1+R$ contains a dicycle of $W$-length $n_1$. Thus $D[V(G_2-y_1+x_1)]$ does not contain a hamiltonian dicycle.  By Corollary 3.5, $D[V(G_2-y_1)]$ is not hamiltonian.  In particular, $y_2y_t\nin A'$. By (24), $d(y_t,G_1)\geq s/2$ as $d(y_1,G_1)=s/2$. If $d(y_t,G_2)\leq (t-1)/2$, then  with $y_t$ in place of $y_1$ in the above argument, we see that $d(y_t,G_1)=s/2$ and so $y_ty_2\in A'$ by (25), \tra. Therefore $d(y_t,G_2)\geq t/2$. With (28), we apply Lemma 2.4 to $y_2y_3\ldots y_t$ in $G_2-y_1$, we see that there exists $a\in\{2,\ldots, t-2\}$ such that for each $i\in\{a+1,a+2,\ldots, t\}$, there exists a hamiltonian path of $G_2$ from $y_1$ to $y_i$ containing $y_1y_2\ldots y_a$, $d(y_i, G_2)\geq t/2$ and $N(y_i,G_2)\bs\{y_a,y_{a+1},\ldots, y_t\}$. Clearly $G[\{y_a,y_{a+1},\ldots, y_t\}]$ is hamiltonian. Thus $a\ne 2$ and so $a\geq 3$. By Lemma 3.6, $y_1y_i\nin E$ for all $i\in\{a+1,a+2,\ldots, t\}$. As $y_ty_2\nin A'$ and by (20), (21), (22) and (27), we obtain that $d(y_1y_t,G_1)=s+1$ and $d(y_1y_t,G_2)=t-1$. It follows that $N_G(y_1,G_2)=\{y_2,y_3,\ldots, y_a\}$. Thus $G_2$ has a hamiltonian path from $y_2$ to $y_t$. By (27) with $y_2$ in place of $y_1$,  we obtain that $y_2y_t\in A'$, \tra. \qed\vs

We are now in a position to complete our proof. By (24) , we may assume w.l.o.g. that $d(y_1,G_1)\geq s/2$. By Lemma 3.6 and Lemma 3.7, $d(y_1,G_2-y_t)\geq t/2$. With (28), applying Lemma 2.4 to $y_{t-1}y_{t-2}\ldots y_1$,  there exists $p\geq \lceil t/2\rceil $ such that for each $i\in\{1,\ldots, p\}$, $N_G(y_i,G_2)\bs\{y_1,\ldots, y_{p+1}\}$, $G_2$ has a hamiltonian path from $y_t$ to $y_i$ containing $y_ty_{t-1}\ldots y_{p+1}$ and  $d(y_i,G_2)\geq d(y_1,G_1)$.   By Lemma 3.6,  $N_G(y_t, G_2)\bs\{y_{p+1},y_{p+1},\ldots, y_{t-1}\}$ and so $d(y_t,G_2)\leq t-1-p\leq \lfloor (t-2)/2\rfloor$. By Lemma 3.7, $d(y_t,G_1)\leq \lfloor (s-1)/2\rfloor$.  Thus $d(y_1,G_1)\geq d(y_1y_t,G_1)-d(y_t,G_1)\geq \lceil (s+1)/2\rceil$ by  (24).  By (27) with $y_i$ in place of $y_1$, $y_ty_i\in A'$ for all $i\in\{1,2,\ldots, p\}$. Thus for all $j\in\{1,2,\ldots, p\}$, $D[V(G_2)-y_j]$ is hamiltonian and so $D[V(G_2)-y_j+x_1]$ is hamiltonian by Corollary 3.5.  We conclude that $D[V(G_1-x_1+y_i+R)]$ does not contain a dicycle of $W$-length $n_1$ for all $i\in\{1,2,\ldots,p\}$. This implies that $a(y_i,x_2x_s)\leq 2$ for all $i\in\{1,2,\ldots, p\}$ and so $a(x_2x_s,y_1y_2\ldots y_p\})\leq 2p$. 

Assume that $x_2\nin W$. By Corollary 3.5, $G_2-y_1+x_1$ is traceable.  By (7), $y_1x_3\nin E$ and $y_1x_s\nin E$. By (7) and Lemma 2.1$(a)$, $d(y_1,P_1-\{x_1,x_2,x_3,x_s\})\leq (s-3)/2$. As $d(y_1,G_1)\geq (s+1)/2$, it follows that  $d(y_1,x_1x_2)=2$ and $d(y_1,P_1-\{x_1,x_2,x_3,x_s\})=(s-3)/2$. Thus $y_1x_s\nin A'$ for otherwise $D[V(G_1-x_1+y_1)]$ is hamiltonian.  By (20), (21), (22) and (27), we obtain that $d(y_1y_t,G_1)= s+1$ and so $d(y_t,G_1)=(s+1)/2$, a contradiction since we have $d(y_t,G_1)\leq \lfloor (s-1)/2\rfloor$ in the above paragraph.  

Therefore $x_2\in W$.  We claim that $a(x_2x_s,R)\leq 2r$  and $x_2x_s\nin A'$. If this is not true, then either $e(z,x_2x_s)\geq 3$ for some $z\in R$ or $x_2x_s\in A'$. Since  $D[V(G_1-x_1+R+y_1)]$ does not contain a dicycle of $W$-length $n_2$, we see that $d(y_1,x_2x_s)\leq 1$ and $d(y_1,x_ix_{i+1})\leq 1$ for all $i\in\{2,\ldots, s-1\}$. By Lemma 2.1$(a)$, this yields that $d(y_1,P_1-x_1)\leq (s-1)/2$. It follows that $d(y_1,G_1)=(s+1)/2$, $y_1x_1\in E$ and $d(y_1,x_2x_s)=1$. Thus either $y_1x_2\nin A'$ or $y_1x_s\nin A'$. By (20), (21), (22) and (27), it follows  that $d(y_1y_t,G_1)= s+1$ and so $d(y_t,G_1)=(s+1)/2$, a contradiction.  Therefore $a(x_2x_s,R)\leq 2r$  and $x_2x_s\nin A'$. Thus $s\geq 4$. 

Next, we claim that $a(x_2x_s,P_1-x_1)\leq 2(s-2)$. If this is not true, then by Lemma 2.3, $D[V(G_1-x_1)]$ has a hamiltonian dicycle $C$. By Lemma 2.2, $d(x_1,C)\leq (s-1)/2$ and  $d(y_1, C)\leq (s-1)/2$. It follows that $d(y_1,C)=(s-1)/2$ and $y_1x_1\in E$.   Since $d(y_1,G_1)=(s+1)/2$, $G_1-x_1+y_1$ is traceable by Lemma 2.1$(b)$.  By Corollary 3.5, $G_2-y_1+x_1$ is traceable. By Lemma 3.2, $(s-1)/2+d(y_1,G_2)\geq d(x_1,G_2)+(s+1)/2-2$. Thus $p\geq d(y_1,G_2)\geq d(x_1,G_2)-1\geq (t+s+r-\la)/2-1$.  Then $d(y_t,G_2)\leq t-1-p\leq (t-s-r+\la)/2$. By Lemma 3.7, $d(y_t,G_1)\leq (s-1)/2$. Thus
\[ 
d(y_t,G)\leq (t-s-r+\la)/2+(s-1)/2+r=(t+r-1+\la)/2<(n-\la)/2.\]
This contradicts (2). 

Therefore $a(x_2x_s,P_1-x_1)\leq 2(s-2)$. Clearly, $a(x_2x_s,x_1)=2$ as $x_1x_s\nin A'$. Let $X=R\cup V(G_1)\cup\{y_1,\ldots,y_p\}$. We obtain that $a(x_2x_t,X)\leq 2r+2(s-2)+2+2p=2(|X|-1)$. By (6), $n-|X|=t-p\geq (t+s+r-\la)/2$. Thus $t\geq 2p+s+r-\la\geq t+s+r-\la>t$, \tra. This proves the theorem.

\end{document}